\documentclass[11pt]{article}
\usepackage{amssymb,url,amsmath,graphicx, amsthm,mathdots}
\usepackage{accents}
\usepackage{tikz}
\usepackage{booktabs}
\usetikzlibrary{positioning, calc, chains}
\usetikzlibrary{shapes.misc}
\usetikzlibrary{shapes.symbols}
\usetikzlibrary{shapes.geometric}
\usetikzlibrary{shapes.arrows}
\usetikzlibrary{fit}
\usetikzlibrary{shadows}

\usepackage{pgfplots}
\usepackage{pgfplotstable}

\usepackage{subcaption}

\usepackage{listings}
\definecolor{DarkBlue}{rgb}{0.00,0.00,0.55}
\definecolor{DarkRed}{rgb}{0.55,0.00,0.00}
\definecolor{DarkGreen}{rgb}{0.00,0.55,0.00}
\definecolor{Bittersweet}{rgb}{1.0, 0.44, 0.37}
\definecolor{Purple}{rgb}{0.5, 0.0, 0.5}

\newcommand{\sech}{\operatorname{sech}}

\newcommand{\citep}[1]{\cite{#1}}
\newcommand{\bfc}{\mathbf{c}}

\newcommand{\bfb}{\mathbf{b}}

\def\hdiv{{H(\mathrm{div})}}

\usepackage{algorithmic}
\usepackage{algorithm}

\title{Irksome: Automating Runge--Kutta time-stepping for finite
  element methods}
\date{}
\author{
  Patrick E.~Farrell\thanks{Mathematical Institute, University
    of Oxford; Woodstock Road, Oxford OX2 6GG, UK. Email:
    patrick.farrell@maths.ox.ac.uk. This work was supported
    by EPSRC grants EP/R029423/1 and EP/V001493/1.}
  \and
    Robert C.~Kirby\thanks{Department of Mathematics, Baylor
    University; One Bear Place \#97328; Waco, TX 76798-7328.
    Email: robert\_kirby@baylor.edu.  This work was supported by
    National Science Foundation grant 1912653.}
\and
 Jorge Marchena-Men\'endez\thanks{Department of Mathematics, Baylor
    University; One Bear Place \#97328; Waco, TX 76798-7328.
    Email: jorge\_marchena1@baylor.edu.}
}
\begin{document}
\maketitle

\begin{abstract}
While implicit Runge--Kutta methods possess high order accuracy and important stability properties, implementation difficulties and the high expense of solving the coupled algebraic system at each time step are frequently cited as impediments.  We present \texttt{Irksome}, a high-level library for manipulating UFL (Unified Form Language) expressions of semidiscrete variational forms to obtain UFL expressions for the coupled Runge--Kutta stage equations at each time step.  \texttt{Irksome} works with the Firedrake package to enable the efficient solution of the resulting coupled algebraic systems.  Numerical examples confirm the efficacy of the software and our solver techniques for various problems.  
\end{abstract}

\section{Introduction}
Many successful high-level finite element packages provide a
domain-specific language that allows users to succinctly specify a
partial differential equation (PDE) in mathematical syntax and obtain
an efficient implementation.  UFL, the Unified Form
Language~\cite{Alnaes:2014}, provides such a Python-based domain-specific
language.  It was originally introduced as part of FEniCS~\cite{Logg:2012} but
has since been adopted by both Firedrake~\cite{Rathgeber:2016} and
DUNE-FEM~\cite{dedner2010generic}.

In addition to enabling the generation of low-level finite element code,
such an abstract interface also provides opportunities for
``outer-loop'' development of algorithms built on top of PDE solves.
For example, \texttt{dolfin-adjoint}~\cite{farrell2013automated} allows users to specify
forward models and cost functionals in UFL and then automatically derives adjoint
and tangent-linear models, key ingredients in data
assimilation, optimization, and sensitivity analysis.
The more recent \texttt{hIPPYlib}~\cite{villa2018hippylib} provides many of
these features but also includes a framework for Bayesian optimization
and certain efficient Hessian approximations.  \texttt{RBniCS}~\cite{ballarin2015rbnics,hesthaven2016certified}
provides several approaches to reduced order modeling, built on top of
FEniCS' interface.

Despite its richness of support for diverse spatial discretizations of various kinds and orders, UFL lacks a comparable abstraction for time-stepping.  Although early
versions of DOLFIN interfaced to multi-adaptive temporal Galerkin
methods~\cite{logg2003multi}, users wishing to solve time-dependent problems
typically write their own time-stepping loops with relatively
elementary methods.  A FEniCS-based library for time abstractions was
developed in~\cite{maddison2014rapid}, although this code mainly aims to
manage the time levels and provide some simple multi-step methods.

In this paper, we present~\texttt{Irksome}, a solution to
the problem of obtaining effective time-discretizations.  Rather than consider all possible kinds of time-stepping algorithms, we restrict ourselves to Runge--Kutta (RK) methods in this project.  Runge--Kutta methods themselves are a vast family covering explicit and implicit methods of various orders of accuracy.   When applications call for some kind of stability or conservation property, there is typically some suitable RK method available.  Moreover, \emph{Butcher tableaux} provide a unified description of RK methods and hence a unified entry point for automating their deployment in conjunction with UFL-based PDE descriptions.
IMEX-type and partitioned Runge--Kutta methods that allow different methods for different terms in the equation can be of use for Hamiltonian systems, but such methods require additional UFL manipulation and we regard them as a future project.  We also expect multi-step, generalized linear methods or other broad families of methods could admit similar deployment.
An alternative approach would be to wrap PETSc's \texttt{TS} package \cite{zhang2018ts}\footnote{This has been recently done in~\url{https://github.com/IvanYashchuk/firedrake-ts}.}, providing access to these other kinds of methods currently available in PETSc.  However, TS currently only provides a small, enumerated set of RK methods, and these do not include multi-stage implicit RK methods.  A major motivation of our work is to study efficient algebraic solvers for the large, coupled algebraic systems such methods yield.

Our work follows the design principle of separating \emph{mechanism} from \emph{policy}~\cite{lampson1976reflections}.  That is, \texttt{Irksome} makes it straightforward to apply a particular RK method to a PDE, but does not comment on which method the user ought to select.  In some ways, this mirrors the agnosticism of UFL and other domain-specific languages for spatial discretizations: UFL provides relatively broad access to classical, mixed, and discontinuous Galerkin methods without editorializing on their relative merits for a particular problem.

At the same time, the discussion of `policy' regarding Runge--Kutta time-stepping of finite element discretizations of PDE is ongoing, and we hope to provide a useful tool in this respect.   Despite their favorable accuracy and stability properties, higher-order fully implicit RK methods have not yet been extensively used in practical applications.  Two reasons for this are typically given.  First, as the algebraic system couples all RK stages, in a traditional finite element code the cost of implementation is quite high; the space and time discretization must be combined in the assembly process.  Second, they have been considered impractically expensive, due to the size of the nonlinear systems to be solved. \texttt{Irksome} addresses both of these points: by automating the application of fully implicit RK methods, the first concern dissolves, and by building on the sophisticated solver infrastructure of Firedrake and PETSc \cite{kirby2018solver}, fast preconditioners can be developed and applied to address the second.

One response to the difficulties of implicit RK methods has been the development of `diagonally implicit' Runge--Kutta methods (DIRKs). These can provide many favorable properties, but with a sequence of smaller algebraic systems for each stage. However, they are necessarily limited to low stage order, and in some cases effective preconditioners for fully implicit methods can be competitive in efficiency with DIRKs.  Mardal, Nilssen, and Staff~\cite{mardal2007order} gave a rigorous analysis of certain block-diagonal preconditioners for fully implicit discretizations of parabolic PDE.  Essentially, the cost of applying the preconditioner for an $s$ stage method is the solution of $s$ independent systems similar to that for a backward Euler step.  When the outer iteration count is small, so is the relative cost of a fully implicit method compared to a DIRK.  Huang \emph{et al}~\cite{huang2019conditioning} combine a diagonal preconditioner for backward Euler with the Jordan form of the Butcher matrix to precondition the higher-order system.  Pazner and Persson~\cite{pazner2017stage} apply an earlier recommendation of Butcher for methods with invertible $A$ matrix~\cite{butcher1976implementation} to the compressible Navier-Stokes equations.  In this technique, they change variables to render the stiffness part of Jacobian block diagonal and hence cheaper to apply and precondition.  They find that with appropriate preconditioning, Radau IIA methods (which are L-stable and have relatively high stage order) can be made more efficient than DIRKs of comparable order.

In  Section~\ref{sec:rk}, we give a brief overview of Runge--Kutta methods.  In particular, we are interested in methods given by a classical Butcher tableau, any of which can be encoded by a few arrays.  In addition to a few general families of collocation methods chosen to highlight different stability features, we have also included a range of classical methods.
In Section~\ref{sec:imp} we describe \texttt{Irksome}, giving an overview of the library itself and some of the internal implementation.  Several examples demonstrating various features are given in Section~\ref{sec:results}, and concluding thoughts are given in Section~\ref{sec:conc}.

\section{Runge--Kutta methods}
\label{sec:rk}
\subsection{Overview}
We first consider ordinary differential equations of the form
\begin{equation}
y^\prime(t) + F(t, y) = 0,
\end{equation}
where $F: (0,T] \times \mathbb{R}^{n} \rightarrow \mathbb{R}^n$,
and the solution $y: (0,T] \rightarrow \mathbb{R}^n$. The equation must also
satisfy some initial condition
\begin{equation}
  y(0) = y_0.
\end{equation}
We will assume that standard conditions for the existence, uniqueness,
and continuous dependence of solutions (e.g.~via the Picard-Lindel\"of
theorem) hold for the discrete ODE.

Given the solution $y(t^n) \equiv y_n $ and some
$t^{n+1} = t^n + \Delta t$,
Runge--Kutta methods approximate $y(t^{n+1})$ by
\begin{equation}
  y^{n+1} = y^n + \Delta t \sum_{i=1}^s b_i k_i,
\end{equation}
where for all $1 \leq i \leq n$ the \emph{stages} $k_i \in \mathbb{R}^n$ satisfy
\begin{equation}
  k_i + F\left(t+c_i\Delta t, y+\Delta t\sum_{i=1}^n A_{ij} k_j
  \right) = 0.
\end{equation}
The numbers $A_{ij}$, $b_i$, and $c_i$ for $1 \leq i, j \leq n$ are in
principle arbitrary but are chosen so that the resulting method has a
given order of accuracy as well other desired properties (e.g.~various
notions of stability or symplecticity).  They are typically organized
in a \emph{Butcher tableau}
\begin{equation}
  \begin{array}{c|c}
    \bfc & A \\ \hline
    & \bfb
  \end{array},
\end{equation}
where the vectors $\bfb,\bfc \in \mathbb{R}^s$ contain the entries $b_i$ and
$c_i$, respectively, and $A \in \mathbb{R}^{s \times s}$ contains the
entries $A_{ij}$.  If $A$ is strictly lower triangular, then the method is explicit -- each stage value can be computed in sequence without recourse to an algebraic system (modulo mass matrices in the variational context). Otherwise, the method is implicit.  In the case of a fully implicit method ($A$ being essentially dense), one must solve an $(ns) \times (ns)$ system of algebraic equations to determine all the stage values simultaneously.  When $A$ is lower triangular but not strictly so (a diagonally implicit method or DIRK), one may solve $s$ consecutive $n \times n$ algebraic systems for the stages.  This is typically more efficient than  solving the single large system for fully implicit methods.

More generally, discretizations of certain PDE (e.g.~mixed formulations of the 
 time-dependent Navier-Stokes) give rise to systems of (explicit) differential-algebraic equations
\begin{equation}
  \begin{split}
    x^\prime(t) + F(t, x, y) & = 0, \\
    G(x, y) & = 0,
  \end{split}
\end{equation}
or even more generally a fully implicit formulation
\begin{equation}
  F(t, x, x^\prime) = 0.
\end{equation}

In this paper, we investigate the application of RK methods to 
finite element spatial discretizations of (generally nonlinear) time-dependent PDE.  For a finite-dimensional function space $V_h$, we consider variational
evolution equations of the form of finding
$u:[0,T) \rightarrow V_h$ such that
  \begin{equation}
    \label{eq:varode}
    (u_t, v) + \mathcal{F}(t, u; v) = 0, 
  \end{equation}
for all $0 < t < T$ and $v \in V_h$.

Here, we assume that $v$ enters into $\mathcal{F}$ linearly but make
no particular assumptions about its dependence on $t$ and $u$, other than
that the discrete system is solvable. We
do not assume that we have a particular spatial discretization -- 
discontinuous Galerkin or other nonconforming methods fit into our
framework as well as standard conforming ones.  

It is useful to consider an even more general problem of finding
$u:[0,T) \rightarrow V_h$ such that
\begin{equation}
\label{eq:generalDAE}
  \mathcal{G}(t, u, u_t; v) = 0,
\end{equation}
for all $0 < t < T$ and $v \in V_h$.  This general formulation allows
some interesting cases (e.g.~Sobolev equations) not covered
in~\eqref{eq:varode}, and working with this abstract form makes our
UFL manipulation more straightforward.

The Runge--Kutta methods we consider have a range of different kinds of stability properties appropriate for various PDE examples we consider later.  We refer the reader to~\cite{hairer2006geometric, wanner1996solving} for an in-depth discussion of these various properties.  The RK methods we consider are $A$-stable, meaning that if applied to the test equation $y^\prime = \lambda y$ with $\mathrm{Re}(\lambda) < 0$, the computed solution tends to zero at infinity for all positive time-steps.  While the famous Dahlquist Barrier Theorem~\cite{dahlquist1963special} says that $A$-stable multi-step methods can be at most second order accurate, $A$-stable RK methods of all orders are known.  $A$-stability potentially allows large time steps to be chosen, provided sufficient accuracy is obtained with them.  For very stiff equations, a stronger notion of stability than $A$-stability is desirable.  So-called $L$-stability requires that the stability function vanish at infinity.
Additionally, some ODE systems possess a monotonicity property that if $y^\prime = f(t, y)$ and $z^\prime = f(t, z)$ are solutions with different initial conditions, then $\| y(t_2) - z(t_2) \| \leq \| y(t_1) - z(t_1) \|$ for $t_1 \leq t_2$.  $B$-stable numerical methods~\cite{wanner1996solving} preserve a discrete analog of this property.

Finally, some ODE possess one or more conserved quantities, and certain RK methods excel in preserving this in the discretization.  For problems with a Hamiltonian structure, \emph{symplectic} methods conserve the Hamiltonian up to a perturbation over exponentially long time scales (the result is stronger for linear problems).  Some methods are additionally known to conserve linear or quadratic invariants, which makes them highly attractive for problems with these features.  We demonstrate such properties for the linear wave equation and the nonlinear Benjamin-Bona-Mahony equations.

\subsection{PDE examples}
To fix ideas, we now consider a range of PDE to serve as motivating
examples.  These examples are chosen to explore two different dimensions of our work.  First, they cover ODE, DAE, and implicit/Sobolev PDE.  Second, these examples call for different kinds of stability properties, highlighting the need for different methods.

\subsubsection{The heat equation}
A model problem for time stepping finite element discretizations is the heat equation 
\begin{equation}
\label{eq:heat}
  \left(u_t, v\right) - \left(\nabla u, \nabla v\right) = \left(f, v \right),
\end{equation}
posed on some domain $\Omega \subset \mathbb{R}^d$ with $d \in \{1,2,3\}$,
together with Dirichlet boundary conditions
$u|_{\partial \Omega} = g(t, \cdot)$.  We let $V_h$ consist of standard continuous
piecewise polynomials defined over a triangulation of $\Omega$.

Backward  Euler is a very common method, highly stable but only first-order accurate, for the heat equation.  Given $u^n$, the approximation to the solution $u$ at time $t_n$, we define $u^{n+1}$ as the solution to the variational problem
\begin{equation}
\label{eq:beheat}
\left( \frac{u^{n+1}-u^n}{\Delta t}, v \right)
+ \left( \nabla u^{n+1}, \nabla v \right) - \left( f\left(t^{n+1}, \cdot \right), v \right) = 0,
\end{equation}
for all $v \in V_h$.  Since $u^{n+1}$ is the unknown value and $u^n$ is data for this problem, we can rearrange this to obtain
\begin{equation}
\label{eq:beheat2}
\left( u^{n+1} , v \right)
+ \Delta t \left( \nabla u^{n+1} , \nabla v \right)
= \left( u^n, v \right) + \Delta t \left( f \left(t^{n+1}, \cdot \right), v \right).
\end{equation}

Although backward Euler is just a Runge--Kutta method with Butcher tableau
\begin{equation}
\label{eq:bebt}
  \begin{array}{c|c}
    1 & 1 \\ \hline
    & 1
  \end{array},
\end{equation}
one typically poses the problem for $u^{n+1}$, as we have done in~\eqref{eq:beheat} or~\eqref{eq:beheat2}, instead of the lone stage $k_1$.   If we apply a generic $s$-stage RK method to~\eqref{eq:heat}, we pose a variational problem for the $s$ stages.  We seek $\{ k_i \}_{i=1}^s \subset V_h$ such that
\begin{equation}
\label{eq:rkheat}
\left( k_i , v_i \right) + \left( \nabla \left( u^n + \Delta t \sum_{j=1}^s a_{ij} k_j \right), \nabla v_i \right) 
- \left(f\left(t + c_i \Delta t , \cdot \right), v_i \right)= 0,
\end{equation}
for each $v_i \in V_h$, and then find $u^{n+1}$ by
\begin{equation}
u^{n+1} = u^n + \sum_{i=1}^s b_i k_i.
\end{equation}

In the case of backward Euler, the variational equation for the stage $k_1$ is similar to~\eqref{eq:beheat}.  If we substitute in~\eqref{eq:bebt}, we find
\begin{equation}
\left( k_1 , v_1 \right) 
+ \left( \nabla \left( u^n + \Delta t k_1 \right), \nabla v_1 \right) 
- \left(f\left(t + \Delta t , \cdot\right), v_1 \right)= 0.
\end{equation}
Now, we can rearrange this to give an equation for $k_1$:
\begin{equation}
\label{eq:bejustk1}
\left( k_1 , v_1 \right) 
+ \Delta t \left( \nabla k_1 , \nabla v_1 \right) 
= \left( f(t + \Delta t, \cdot), v_1 \right) - \left( \nabla u^n , \nabla v_1 \right).
\end{equation}
Although the right-hand side is different, the quantity $k_1$ here can be put into $1-1$ correspondence with $u^{n+1}$ from~\eqref{eq:beheat} by 
$k_1 \leftrightarrow \tfrac{u^{n+1} - u^n}{\Delta t}$.  Importantly, we note the equality of the bilinear form defining $u^{n+1}$ in~\eqref{eq:beheat2} to that defining $k_1$ in~\eqref{eq:bejustk1}.

While the Runge--Kutta formulation~\eqref{eq:rkheat} is quite a bit different than~\eqref{eq:beheat}, we can programmatically obtain it from~\eqref{eq:heat} by, for each stage, substituting $k_i$ in for $\tfrac{\partial u_i}{\partial t}$, $u_n + \Delta t \sum_{j=1}^s a_{ij} k_j$ for $u$, $v_i$ for $v$, and $t + c_i \Delta t$ for t.  This simple insight forms the basis of the implementation, which we discuss in section \ref{sec:imp}.

Strongly-enforced boundary conditions for the semidiscrete problem must be converted to boundary conditions for each Runge--Kutta stage.  The stage $k_i$ approximates the time derivative $u_t$ at time $t^n + c_i \Delta t$, so if we pose a boundary condition of the form
\[
u|_\Gamma = g(t, \cdot),
\]
for some boundary segment $\Gamma \subset \partial \Omega$, we impose the corresponding stage boundary condition
\begin{equation}
(k_i)|_\Gamma = \tfrac{\partial g}{\partial t}(t^n+c_i\Delta t, \cdot).
\end{equation}
We note that weakly-enforced boundary conditions (such as Neumann and Robin for standard Galerkin methods) require no special treatment, just the evaluation of any time-dependent data at the correct stage times.

The discrete heat equation is a prototypical stiff system.  Explicit time-stepping methods require a time step with $\Delta t = \mathcal{O}(h^2)$, where $h$ is the mesh size.  $A$-stability allows us to put $\Delta t = \mathcal{O}(h)$, with a constant possibly larger than 1.  Also, because the stiffness increases under mesh refinement, $L$-stability can be attractive.  Hence, we would expect Gauss-Legendre methods to be passable but perhaps to obtain better results from RadauIIA and LobattoIIIC families for a given order.  While the heat equation does possess a contractive property, $B$-stability does not seem as important in practice for these problems.

Generalizing the class of problems we consider, the mixed form of the heat equation gives rise to a system of differential algebraic equations (DAE).  Introducing the new variable $\sigma = -\nabla u$ in~\eqref{eq:heat}, we obtain the system of PDE
\begin{equation}
\begin{split}
u_t + \nabla \cdot \sigma - f & = 0, \\
\sigma + \nabla u & = 0.
\end{split}
\end{equation}
To discretize this, we take $W^1_h \subset L^2$ as the space of discontinuous piecewise polynomials of degree $k$ over a triangulation of $\Omega$ and $W^2_h \subset \hdiv$ a suitable mixed approximating space and define $V_h = W^1_h \times W^2_h$.  We then seek $(u, \sigma): (0, T] \rightarrow V_h$ such that
\begin{equation}
\label{eq:mixedheat}
\begin{split}
\left( u_t, v \right) + \left( \nabla \cdot \sigma, v \right)- \left(f, v\right) & = 0, \\
\left( \sigma, w \right) - \left(u, \nabla \cdot w \right) & = 0,
\end{split}
\end{equation}
for all $(v, w) \in V_h$.  Since only $u_t$ appears in the equation and not $\sigma_t$, we have a differential-algebraic system rather than just a system of ODE.  Still, we can apply a generic RK method to~\eqref{eq:mixedheat} to obtain a coupled system of variational problems for all $s$ stages.  For each $1 \leq i \leq s$, we seek $(k^u_i, k^\sigma_i) \in V_h$ such that
\begin{equation}
\label{eq:rkmixedheat}
\begin{split}
\left( k^u_i, v_i \right) 
+ \left( \nabla \cdot \left( \sigma^n + \Delta t \sum_{j=1}^s a_{ij} k^\sigma_{j} \right), v_i \right)- \left(f(t_n + c_i \Delta t, \cdot), v_i\right) & = 0, \\
\left( \left( \sigma^n + \Delta t \sum_{j=1}^s a_{ij} k^\sigma_{j} \right), w_i \right) 
- \left(\left(u^n + \Delta t \sum_{j=1}^s a_{ij} k^u_{j} \right), \nabla \cdot w_i \right) & = 0,
\end{split}
\end{equation}
for all $(v_i, w_i) \in V_h$.  %This is a DAE, although in a rather explicit formulation.

In some sense, DAE are ``infinitely stiff'', and so we expect $A$- and $L$-stability to be quite important for these problems. We therefore expect RadauIIA or LobattoIIIC methods to outperform Gauss-Legendre.

\subsubsection{The wave equation}
We can write the wave equation $u_{tt} - \Delta u = 0$ as a
first-order variational system
\begin{equation}
\label{eq:mixedwave}
  \begin{split}
    (u_t, v) + (\nabla \cdot \sigma, v) & = 0, \\
    (\sigma_t, w) - (u, \nabla \cdot w) & = 0.
  \end{split}
\end{equation}
Following~\cite{geveci, kirbykieu}, the flux variable $\sigma \in
H(\mathrm{div})$ and $u \in L^2$.  As above, we take $W^1_h$ as the space of discontinuous piecewise polynomials of degree $k$ over a triangulation of $\Omega$ and $W^2_h \subset \hdiv$ a suitable mixed approximating space and define $V_h = W^1_h \times W^2_h$.

The resulting ODE tend to be far less stiff than those for
heat equation.  Energy conservation rather than stiffness tends to be the central issue for time stepping.  Standard integrators like forward or backward Euler dramatically fail to conserve energy, and one typically requires some kind of symplectic integrator.  In~\cite{kirbykieu}, Kirby analyzed a first-order symplectic Euler time-stepping scheme for this problem, and Kernell \&Kirby~\cite{kernell2020preconditioning} studied preconditioners for a Crank-Nicolson time discretization of a slightly more general equation.

Applying a generic Runge--Kutta method to~\eqref{eq:mixedwave} gives 
\begin{equation}
\label{eq:rkmixedwave}
\begin{split}
\left( k^u_i, v_i \right) 
+ \left( \nabla \cdot \left( \sigma^n + \Delta t \sum_{j=1}^s a_{ij} k^\sigma_{j} \right), v_i \right)- \left(f(t_n + c_i \Delta t, \cdot), v_i\right) & = 0, \\
\left( k^\sigma_i, w_i \right)
- \left(\left(u^n + \Delta t \sum_{j=1}^s a_{ij} k^u_{j} \right), \nabla \cdot w_i \right) & = 0.
\end{split}
\end{equation}
Some Runge--Kutta families, such as Gauss-Legendre methods, provide high-order, $A$-stable, and symplectic methods that preserve system energy for linear problems (and nearly so for nonlinear ones)~\cite{hairer2006geometric}.  The damping features that make RadauIIA and LobattoIIIC quite suitable for the heat equation turn out to be a major weakness for the wave equations.  Despite the different RK method to be preferred, we note that the discrete system~\eqref{eq:rkmixedwave} has quite a similar structure to~\eqref{eq:rkmixedheat}.

\subsubsection{Nonlinear examples}
As an example of a PDE system giving rise to a nonlinear DAE system, we consider  the incompressible Navier-Stokes equations
\begin{equation}
  \label{eq:nse}
  \begin{split}
    u_t - \nu \Delta u + u \cdot \nabla u + \nabla p & = 0 ,\\
    \nabla \cdot u & = 0,
  \end{split}
\end{equation}
where $u$ is the vector-valued fluid velocity and $p$ the pressure, and the parameter $\nu$ is the kinematic viscosity.  Among the many issues these equations present, one must choose suitably compatible spaces for velocity and pressure to enforce the divergence constraint.

We also consider the Benjamin-Bona-Mahoney~\cite{benjamin1972model} (BBM) equation
\begin{equation}
  \label{eq:bbm}
 u_t + u_x + u u_x - u_{txx} = 0,
\end{equation}
typically posed on $\Omega = \mathbb{R}$ or a periodic domain.
It has multivariate extensions and is also similar to more complex models arising in magma dynamics~\cite{simpson2011solitary}.  For our purposes, it is a nonlinear Sobolev-type equation, with spatial derivatives acting on time derivatives in the $u_{txx}$ term.

The BBM equation has a Hamiltonian structure and three polynomial invariants:
\begin{eqnarray}
  \label{eq:I}
  I_1 = & \int u \, dx, \\
  I_2 = & \int u^2 + \left( u_x \right)^2 \, dx, \\
  I_3 = & \int \left( u_x \right)^2 + \tfrac{1}{3} u^3 \, dx.
\end{eqnarray}
Each of these quantities remains constant over time.  Note that $I_1$ is a linear invariant, $I_2$ quadratic, and $I_3$ cubic.  The BBM equation also supports solitary (but non-soliton) wave solutions.

%% An example that includes both nonlinearity and more extreme
%% stiffness is the fourth-order Cahn--Hilliard model of phase
%% separation:
%% \begin{equation}
%%   \frac{\partial c}{\partial t} - \nabla \cdot M \left(\nabla\left(\frac{d f}{d c}
%%           - \lambda \nabla^{2}c\right)\right) = 0 \quad {\rm in}
%%           \ \Omega.
%% \end{equation}
%% Here $f$ is some typically non-convex function (in our case, we take $f(c)
%% = 100c^2(1-c)^2$), and $\lambda$ and $M$ are
%% scalar parameters controlling rates.  Although $M=M(c)$ (the so-called
%% degenerate mobility case) is possible, we have considered just
%% constant $M$ in our examples.
%% The system is closed with the boundary conditions
%% \begin{align}
%% M\left(\nabla\left(\frac{d f}{d c} - \lambda \nabla^{2}c\right)\right) \cdot n &= 0 \quad {\rm on} \ \partial\Omega, \\
%% M \lambda \nabla c \cdot n &= 0 \quad {\rm on} \ \partial\Omega.
%% \end{align}
%% Although this equation is often re-rewritten as a system of two
%% second-order equations~\citep{barrett1999finite}, our work enabling $C^1$-conforming elements in Firedrake~\citep{finat-zany} allows us to tackle the primal formulation.

\subsection{Algebraic systems}
A long-standing critique of higher-order implicit RK methods, and argument for DIRKs, is the size and complexity of the algebraic systems required to be solved at each time step.  We now discuss some of these issues and attempt to address them in the context of the heat equation with a 3-stage method.  We define $M$ and $K$ to be the standard finite element mass and stiffness matrices, with
\begin{equation}
M_{ij} = \int_{\Omega} \psi_i \psi_j \, dx, \ \ \
K_{ij} = \int_{\Omega} \nabla \psi_i \cdot \nabla \psi_j \, dx
\end{equation}
where $\{\psi_i\}_{i=1}^{\dim V_h}$ is a finite element basis.  The variational problem~\eqref{eq:rkheat} gives rise to a block algebraic system of the form
\begin{equation}
\label{eq:blockmat}
\left(
\begin{bmatrix}
M & 0 & 0 \\
0 & M & 0 \\
0 & 0 & M
\end{bmatrix}
+ \Delta t
\begin{bmatrix}
a_{11} K & a_{12} K & a_{13} K \\
a_{21} K & a_{22} K & a_{23} K \\
a_{31} K & a_{32} K & a_{33} K
\end{bmatrix}
\right)
\begin{bmatrix}
\mathrm{k}_1\\ \mathrm{k}_2 \\ \mathrm{k}_3
\end{bmatrix}
= \begin{bmatrix} \mathrm{f}_1 \\ \mathrm{f}_2 \\ \mathrm{f}_3 \end{bmatrix},
\end{equation}
or equivalently (and true for general $s$-stage methods)
\begin{equation}
  \label{eq:blockmatgen}
\left( I \otimes M + \Delta t A \otimes K \right) \mathrm{\bf k} = \mathrm{\bf f}.
\end{equation}

Following~\cite{mardal2007order}, the block diagonal of this system makes an excellent preconditioner, at least in the case of parabolic problems such as the heat equation:
\begin{equation}
\label{eq:P}
P = 
\begin{bmatrix}
M + a_{11} \Delta t K& 0 & 0 \\
0 & M + a_{22} \Delta t K & 0 \\
0 & 0 & M + a_{33} \Delta t K
\end{bmatrix}.
\end{equation}

The Jacobian for nonlinear problems will similarly couple together all of the RK stages, although without the Kronecker product structure in the stiffness matrix.  We note that if a DIRK method is used so that $A$ is lower-triangular,~\eqref{eq:blockmat} becomes
\begin{equation}
\label{eq:DIRKblockmat}
\left(
\begin{bmatrix}
M & 0 & 0 \\
0 & M & 0 \\
0 & 0 & M
\end{bmatrix}
+ \Delta t
\begin{bmatrix}
a_{11} K & 0 & 0 \\
a_{21} K & a_{22} K & 0 \\
a_{31} K & a_{32} K & a_{33} K
\end{bmatrix}
\right)
\begin{bmatrix}
\mathrm{k}_1\\ \mathrm{k}_2 \\ \mathrm{k}_3
\end{bmatrix}
= \begin{bmatrix} \mathrm{f}_1 \\ \mathrm{f}_2 \\ \mathrm{f}_3 \end{bmatrix}.
\end{equation}
In this case, the triangular structure means one can proceed by forward substitution, solving in turn for each of the stages.
Moreover, an \emph{explicit} method (typically not recommended for the heat equation) would have $A$ strictly lower triangular and the linear system
\begin{equation}
\label{eq:ERKblockmat}
\left(
\begin{bmatrix}
M & 0 & 0 \\
0 & M & 0 \\
0 & 0 & M
\end{bmatrix}
+ \Delta t
\begin{bmatrix}
0 & 0 & 0 \\
a_{21} K & 0 & 0 \\
a_{31} K & a_{32} K & 0
\end{bmatrix}
\right)
\begin{bmatrix}
\mathrm{k}_1\\ \mathrm{k}_2 \\ \mathrm{k}_3
\end{bmatrix}
= \begin{bmatrix} \mathrm{f}_1 \\ \mathrm{f}_2 \\ \mathrm{f}_3 \end{bmatrix},
\end{equation}
so that we only must invert the mass matrix once per stage (and do mat-vec and axpy-type operations) to obtain all of the stages. 

We remark that the Firedrake solver infrastructure~\cite{kirby2018solver} makes it possible to apply the (rather large) operator in a matrix-free fashion, saving large amounts of memory, and apply some kind of multigrid strategy to the diagonal blocks if the preconditioner~\eqref{eq:P} is used.  Moreover, one can apply a block lower-triangular preconditioner (via PETSc's \texttt{FIELDSPLIT}~\cite{brown2012composable}) for DIRK and explicit methods, which exactly solves the linear systems in a single iteration, provided that the diagonal blocks are accurately inverted.  Hence, the range of RK methods enabled in Irksome composes with the existing Firedrake code stack to enable efficient implementations.

Rather than splitting the stages per~\cite{mardal2007order}, it should also be possible to apply an unsplit or monolithic multigrid scheme using some kind of blockwise smoothing.  The theory of such methods is quite undeveloped, but this could be an interesting future direction.

\section{Irksome}
\label{sec:imp}
Irksome is a Python library for manipulating UFL for semidiscrete variational forms and strongly-enforced boundary conditions into UFL for fully discrete Runge--Kutta methods.  It is available on GitHub under the Firedrake project umbrella at \url{https://github.com/firedrakeproject/Irksome}, and can be installed as part of Firedrake by adding the \lstinline{--install irksome} option to the Firedrake installation script.  

Because exploiting the structure of algebraic system is a key ingredient for obtaining efficient solvers, it is crucial that we maintain this symbolic information in implicit Runge--Kutta methods.  Consequently, Irksome's approach based on UFL manipulation seems to afford richer opportunities than a black-box RK library could. As concrete examples, maintaining the symbolic structure enables the exploration of monolithic multigrid and advanced field-split (Schur complement) preconditioners to the multi-stage system, and allows for efficient vectorized matrix-free implementations.

Irksome's main actions are performed by a function \texttt{getForm} that takes UFL for the time-dependent variational form, a Butcher tableau, the current time and time step (both stored as UFL \texttt{Constant} objects), and a UFL \texttt{Coefficient}.  Irksome also provides a collection of Runge--Kutta methods in a separate module (see Subsection~\ref{ssec:RKmethods}).  The function \texttt{getForm} returns UFL for the coupled, multi-stage method and boundary conditions.
We also provide a convenience class, \texttt{TimeStepper} that handles interactions with \texttt{getForm}, boundary conditions, and variational solvers and advances the solution forward in time.  This class and the module containing Butcher tableaux are the main user entry points for Irksome.

\subsection{UFL manipulation}
To illustrate what Irksome automates, we think of the semidiscrete homogeneous heat equation, which can be written in UFL (with our support of a \texttt{Dt} operator) as
\begin{lstlisting}
F = inner(Dt(u), v) * dx + inner(grad(u), grad(v)) * dx
\end{lstlisting}
Typically, we would express backward Euler for the heat equation in UFL as
\begin{lstlisting}
F = inner((unew - u)/dt, v) * dx + inner(grad(unew), grad(v)) * dx
\end{lstlisting}
However, Runge--Kutta methods require specifying the variational form the update stages satisfy rather than the value at the next time step.  To fix ideas, consider the two-stage LobattoIIIC method given by the Butcher tableau
\begin{equation}
  \begin{array}{c|cc}
    0 & {1}/{2} & -{1}/{2}  \\
    1 & {1}/{2} & {1}/{2} \\ \hline
      & {1}/{2} & {1}/{2}
  \end{array}.
\end{equation}
Following~\eqref{eq:rkheat}, we could write the variational form for the two stages as in Listing~\ref{fig:rkufl}, and solve the variational problem \lstinline{F==0} for the unknown \lstinline{k}, and use it to update the solution.

\begin{figure}
\begin{lstlisting}[captionpos=b, language=Python, caption={Sample UFL for a two-stage Runge--Kutta discretization of the heat equation.}, label={fig:rkufl}]
# u is a given Function in V containing solution at time n
# dt is a Constant holding the current time step value
Vbig = V * V
k = Function(Vbig)
k0, k1 = split(k)
v0, v1 = TestFunctions(Vbig)
u0 = u + Constant(0.5) * dt * k0 + Constant(-0.5) * dt * k1
u1 = u + Constant(0.5) * dt * k0 + Constant(0.5) * dt * k1

F = (inner(k0, v0) * dx + inner(k1, v1) * dx +
     inner(grad(u0), grad(v0)) * dx + inner(grad(u1), grad(v1)) * dx)
\end{lstlisting}
\end{figure}

Using a different Butcher tableau leads to a similar variational problem, simply replacing the entry of the entries of $A$.  More generally, we also must utilize $c$ to evaluate explicitly time-dependent expressions (say, in material properties or forcing terms).  Because the transformation from the semidiscrete form of $F$ to its Runge--Kutta variational form is quite mechanical, it can be mechanized.

Given the UFL expression for $\mathcal{G}$ in~\eqref{eq:generalDAE}, we loop over RK stages \lstinline{i} and use UFL's symbolic manipulation facilities to replace time derivatives of each unknown fields with the RK stage variable $k_i$.  Then unknown fields (without time derivatives applied) $u$ are replaced with $u+(\Delta t)\sum_{j=1}^s A_{ij} k_j$, and the test function with the test function $v_i$.  Explicit appearances of the time $t$ are replaced with $t+c_i \Delta t$.  The resulting substitutions are summed over stages $i$ to produce a monolithic variational form much like in Listing~\ref{fig:rkufl}.  Similarly, substitution is performed on boundary condition values.  Some care must be taken in both the variational form and boundary conditions to handle mixed problems where the underlying function space is the Cartesian product of multiple spaces.  Also, in our current implementation, we assume that time derivatives are applied only to unknown fields and not to algebraic combinations of expressions and fields (e.g. one must write \lstinline{2*t*u + t**2 * Dt(u)} rather than \lstinline{Dt(t**2 u)} and \lstinline{2*Dt(u)*u} rather than \lstinline{Dt(u**2)}).  This does not seem to be a major limitation, but could be relaxed by first applying a transformation expanding such derivatives before applying substitutions.

\subsection{Currently available RK methods}
\label{ssec:RKmethods}
The universe of known Runge--Kutta methods is vast.  While any RK method can be implemented by passing the appropriate Butcher tableau, Irksome provides implementations of many familiar classical methods, with a particular focus on fully implicit collocation-type methods. These typically have an adjustable parameter for the order, much like the function space constructors in FEniCS and Firedrake.
 
Among collocation-type Runge--Kutta methods, we have implemented Gauss-Legendre (which includes implicit midpoint), LobattoIIIA (which includes Crank-Nicolson), and RadauIIA (which includes backward Euler) methods.  We use FIAT~\cite{Kirby:2004} to obtain the interpolating nodes for $c$, and then explicit formulae (evaluated by numerical integration) give the entries of $A$ and $b$ in the Butcher tableau.  Since FIAT can compute the quadrature points to arbitrary degree, we similarly have arbitrary-order implementation of these methods.  These methods are all $A$-stable.  The Gauss-Legendre methods are $B$-stable and symplectic, but not $L$-stable, while the RadauIIA methods are $L$-stable.
We have also implemented general-order LobattoIIIC methods, which have the same $b$ and $c$ arrays as LobattoIIIA but a different $A$ matrix.  They are not quite standard collocation methods, but are $L$-stable and $B$-stable as well as $A$-stable and frequently recommended for stiff problems.

In addition to these methods, we provide Butcher tableaux for a range of other classical methods (forward Euler, an SSP method~\cite{gottlieb2001strong}, explicit midpoint/trapezoid rules, the classical fourth-order method, and some L-stable DIRKs~\cite{alexander1977diagonally}).

\section{Examples and numerical results}
\label{sec:results}
\subsection{The heat equation}
In this section, we confirm the accuracy of our generated methods for the heat equation via the method of manufactured solutions.  We  let $\Omega=[0,1]\times[0,1]$ be the unit square and select the forcing function $f$ and Dirichlet boundary conditions such that the true solution is $u(x, y, t) = e^{-t} \sin(\pi x) \cos(\pi y)$.

For the primal form of this equation, we divide $\Omega$ into an $N \times N$ mesh of squares for $N=8,16,32,64,128$ and take $V_h$ to be the space of cubic serendipity finite elements.  On smooth enough solutions, these lead to fourth order $L^2$ error and third order $H^1$ error, so we focus on time-stepping methods of order at least three -- Gauss-Legendre(2) and LobattoIIIC(3) (both of order 4) and two- and three-stage RadauIIA methods of order 3 and 5, respectively.

\begin{figure}[h]
    \begin{subfigure}[l]{0.475\textwidth}
      \caption{Gauss-Legendre(2)}
      \begin{tikzpicture}[scale=0.88]
      \begin{loglogaxis}[xlabel={$N$}, ylabel={$\|u-u_h\|/\|u\|$},
             ylabel near ticks, ymax=1.e-2, ymin=1.e-12,
             legend pos=south west, legend style={font=\tiny} ,
             cycle list name=color list]
        \addplot[dashed, red]
        table [x=N,y=l2, col sep=comma]{primal_heat_S_3_GL2_1.csv};
        \addlegendentry{$L^2$, $\Delta t=1/N$}
        \addplot[densely dashed, blue] 
        table [x=N,y=l2, col sep=comma]{primal_heat_S_3_GL2_4.csv};
        \addlegendentry{$L^2$, $\Delta t=4/N$}
        \addplot[loosely dashed, green]
        table [x=N,y=l2, col sep=comma]{primal_heat_S_3_GL2_16.csv};
        \addlegendentry{$L^2$, $\Delta t=16/N$}
        \addplot[dotted, purple]
        table [x=N,y=h1, col sep=comma]{primal_heat_S_3_GL2_1.csv};
        \addlegendentry{$H^1$, $\Delta t=1/N$}
        \addplot[loosely dotted, orange] 
        table [x=N,y=h1, col sep=comma]{primal_heat_S_3_GL2_4.csv};
        \addlegendentry{$H^1$, $\Delta t=4/N$}
        \addplot[densely dotted, pink]
        table [x=N,y=h1, col sep=comma]{primal_heat_S_3_GL2_16.csv};
        \addlegendentry{$H^1$, $\Delta t=16/N$}
        \addplot [domain=32:128] {3/pow(x,3)} node[above, midway, yshift=-1pt, anchor=south west] {$\mathcal{O}(h^3)$};
        \addplot [domain=32:128] {0.002/pow(x,4)} node[above, midway, yshift=-1pt, anchor=south west] {$\mathcal{O}(h^{4})$};
      \end{loglogaxis}
    \end{tikzpicture}
    \label{subfig:gl2}
    \end{subfigure}
    \hspace{0.04\textwidth}
    \begin{subfigure}[l]{0.475\textwidth}
      \caption{LobatoIIIC(3)}
    \begin{tikzpicture}[scale=0.88]
      \begin{loglogaxis}[xlabel={$N$}, ylabel={$\|u-u_h\|/\|u\|$},
             ylabel near ticks, ymax=1.e-2, ymin=1.e-12,
             legend pos=south west, legend style={font=\tiny} ,
             cycle list name=color list]
        \addplot[dashed, red]
        table [x=N,y=l2, col sep=comma]{primal_heat_S_3_L3_1.csv};
        \addlegendentry{$L^2$, $\Delta t=1/N$}
        \addplot[densely dashed, blue] 
        table [x=N,y=l2, col sep=comma]{primal_heat_S_3_L3_4.csv};
        \addlegendentry{$L^2$, $\Delta t=4/N$}
        \addplot[loosely dashed, green]
        table [x=N,y=l2, col sep=comma]{primal_heat_S_3_L3_16.csv};
        \addlegendentry{$L^2$, $\Delta t=16/N$}
        \addplot[dotted, purple]
        table [x=N,y=h1, col sep=comma]{primal_heat_S_3_L3_1.csv};
        \addlegendentry{$H^1$, $\Delta t=1/N$}
        \addplot[loosely dotted, orange] 
        table [x=N,y=h1, col sep=comma]{primal_heat_S_3_L3_4.csv};
        \addlegendentry{$H^1$, $\Delta t=4/N$}
        \addplot[densely dotted, pink]
        table [x=N,y=h1, col sep=comma]{primal_heat_S_3_L3_16.csv};
        \addlegendentry{$H^1$, $\Delta t=16/N$}
        \addplot [domain=32:128] {2.0/pow(x,3)} node[above, midway, yshift=-1pt, anchor=south west] {$\mathcal{O}(h^3)$};
        \addplot [domain=32:128] {0.002/pow(x,4)} node[above, midway, yshift=-1pt, anchor=south west] {$\mathcal{O}(h^{4})$};
      \end{loglogaxis}
    \end{tikzpicture}
    \label{subfig:l3}
    \end{subfigure}
    \\
    \begin{subfigure}[l]{0.475\textwidth}
    \caption{RadauIIA(2)}
    \begin{tikzpicture}[scale=0.88]
      \begin{loglogaxis}[xlabel={$N$}, ylabel={$\|u-u_h\|/\|u\|$},
             ylabel near ticks, ymax=1.e-2, ymin=1.e-12,
             legend pos=south west, legend style={font=\tiny} ,
             cycle list name=color list]
        \addplot[dashed, red]
        table [x=N,y=l2, col sep=comma]{primal_heat_S_3_R2_1.csv};
        \addlegendentry{$L^2$, $\Delta t=1/N$}
        \addplot[densely dashed, blue] 
        table [x=N,y=l2, col sep=comma]{primal_heat_S_3_R2_4.csv};
        \addlegendentry{$L^2$, $\Delta t=4/N$}
        \addplot[loosely dashed, green]
        table [x=N,y=l2, col sep=comma]{primal_heat_S_3_R2_16.csv};
        \addlegendentry{$L^2$, $\Delta t=16/N$}
        \addplot[dotted, purple]
        table [x=N,y=h1, col sep=comma]{primal_heat_S_3_R2_1.csv};
        \addlegendentry{$H^1$, $\Delta t=1/N$}
        \addplot[loosely dotted, orange] 
        table [x=N,y=h1, col sep=comma]{primal_heat_S_3_R2_4.csv};
        \addlegendentry{$H^1$, $\Delta t=4/N$}
        \addplot[densely dotted, pink]
        table [x=N,y=h1, col sep=comma]{primal_heat_S_3_R2_16.csv};
        \addlegendentry{$H^1$, $\Delta t=16/N$}
        \addplot [domain=32:128] {2.0/pow(x,3)} node[above, midway, yshift=-1pt, anchor=south west] {$\mathcal{O}(h^3)$};
        \addplot [domain=32:128] {0.002/pow(x,4)} node[above, midway, yshift=-1pt, anchor=south west] {$\mathcal{O}(h^{4})$};
      \end{loglogaxis}
    \end{tikzpicture}
    \label{subfig:r2}
    \end{subfigure}
    \hspace{0.04\textwidth}
    \begin{subfigure}[l]{0.475\textwidth}
    \caption{RadauIIA(3)}
    \begin{tikzpicture}[scale=0.88]
      \begin{loglogaxis}[xlabel={$N$}, ylabel={$\|u-u_h\|/\|u\|$},
             ylabel near ticks, ymax=1.e-2, ymin=1.e-12,
             legend pos=south west, legend style={font=\tiny} ,
             cycle list name=color list]
        \addplot[dashed, red]
        table [x=N,y=l2, col sep=comma]{primal_heat_S_3_R3_1.csv};
        \addlegendentry{$L^2$, $\Delta t=1/N$}
        \addplot[densely dashed, blue] 
        table [x=N,y=l2, col sep=comma]{primal_heat_S_3_R3_4.csv};
        \addlegendentry{$L^2$, $\Delta t=4/N$}
        \addplot[loosely dashed, green]
        table [x=N,y=l2, col sep=comma]{primal_heat_S_3_R3_16.csv};
        \addlegendentry{$L^2$, $\Delta t=16/N$}
        \addplot[dotted, purple]
        table [x=N,y=h1, col sep=comma]{primal_heat_S_3_R3_1.csv};
        \addlegendentry{$H^1$, $\Delta t=1/N$}
        \addplot[loosely dotted, orange] 
        table [x=N,y=h1, col sep=comma]{primal_heat_S_3_R3_4.csv};
        \addlegendentry{$H^1$, $\Delta t=4/N$}
        \addplot[densely dotted, pink]
        table [x=N,y=h1, col sep=comma]{primal_heat_S_3_R3_16.csv};
        \addlegendentry{$H^1$, $\Delta t=16/N$}
        \addplot [domain=32:128] {2/pow(x,3)} node[above, midway, yshift=-1pt, anchor=south west] {$\mathcal{O}(h^3)$};
        \addplot [domain=32:128] {0.002/pow(x,4)} node[above, midway, yshift=-1pt, anchor=south west] {$\mathcal{O}(h^{4})$};
      \end{loglogaxis}
    \end{tikzpicture}
    \label{subfig:r3}
    \end{subfigure}
    \label{fig:primalheat}
    \caption{Relative $L^2$ and $H^1$ errors at $t=2$ on an $N \times N$ mesh using various time stepping schemes and CFL numbers.}
\end{figure}
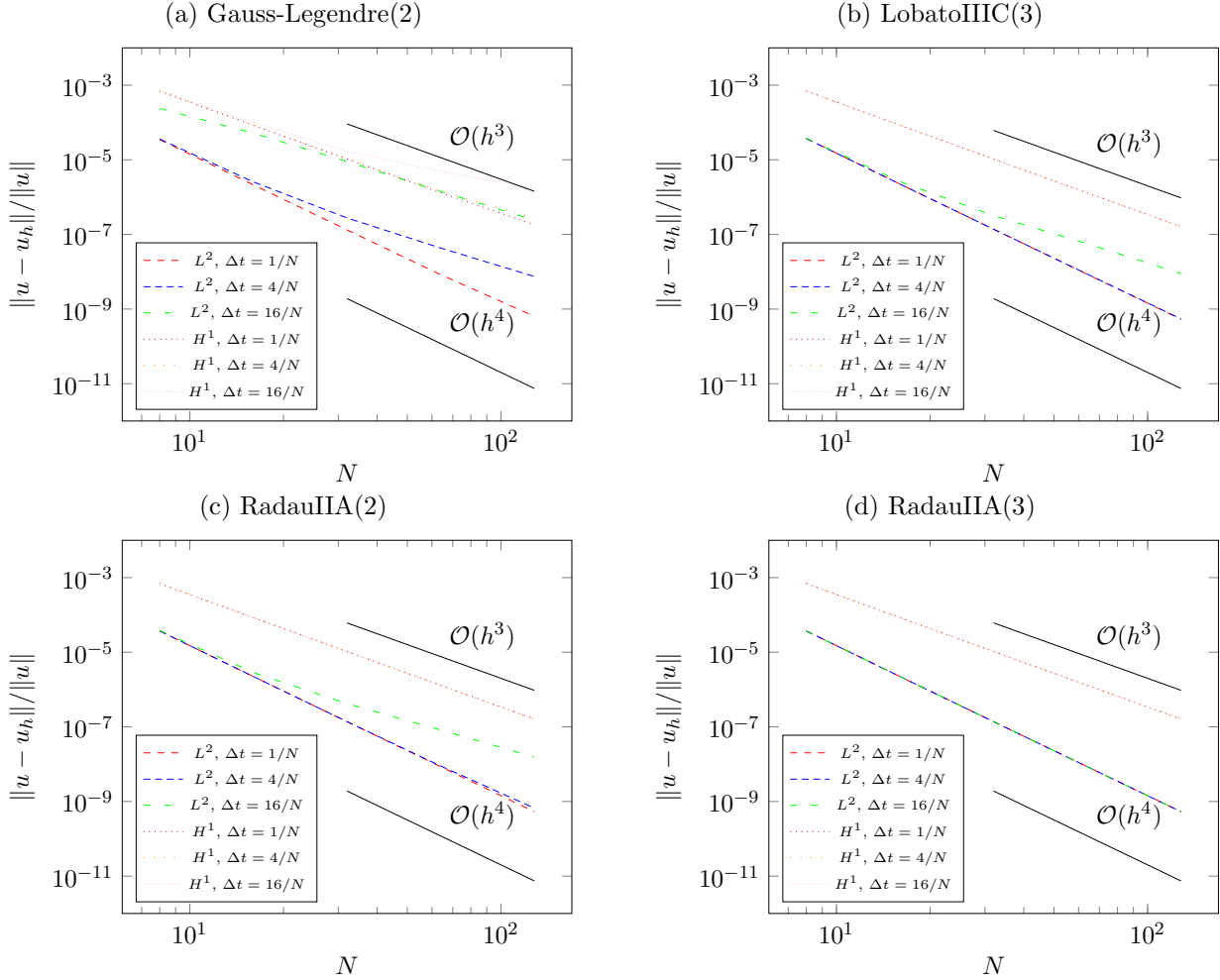

The numerical results illustrate the effect of \emph{order reduction}~\cite{wanner1996solving}.  For stiff problems, collocation-type methods do not yield the full nominal order of accuracy (e.g.~order $2s$ for an $s$-stage Gauss-Legendre method) but reduce to the maximum truncation error per-stage (e.g.~$s$ for an $s$-stage Gauss-Legendre method).  For time steps large enough for temporal error to be significant relative to spatial error, we see a reduction in the observed accuracy.  For example, Figure~\ref{subfig:gl2} shows fourth order accuracy in $L^2$ using the two-stage Gauss-Legendre method when the time step is very small, but this degenerates to stiff order of two.  A similar reduction of the $H^1$ convergence rate below third order is observed for larger time steps.  The $L^2$ error using LobattoIIIC(3) also reduces for larger time steps, but the $H^1$ error stays at order three (which is both the stage order and the spatial order in this norm).  The RadauIIA(2) method with stiff order two has order reduction in $L^2$ norm, but none is observed in $H^1$ (presumably the time truncation error is small enough not to affect the overall convergence order).  The RadauIIA(3) method has stage order three and so does not affect the $H^1$ error, and its time truncation error is small enough that, for the time steps we consider, we observe full accuracy in $L^2$ as well.

\begin{figure}[h]
    \begin{subfigure}[l]{0.475\textwidth}
      \caption{Gauss-Legendre(2)}
      \begin{tikzpicture}[scale=0.88]
      \begin{loglogaxis}[xlabel={$N$}, ylabel={$\|u-u_h\|/\|u\|$},
             ylabel near ticks, ymax=1.e-2, ymin=1.e-9,
             legend pos=south west, legend style={font=\tiny} ,
             cycle list name=color list]
        \addplot[dashed, red]
        table [x=N,y=l2u, col sep=comma]{mixed_heat_simplex_3_GL2_1.csv};
        \addlegendentry{$L^2$, $\Delta t=1/N$}
        \addplot[densely dashed, blue] 
        table [x=N,y=l2u, col sep=comma]{mixed_heat_simplex_3_GL2_4.csv};
        \addlegendentry{$L^2$, $\Delta t=4/N$}
        \addplot[loosely dashed, green]
        table [x=N,y=l2u, col sep=comma]{mixed_heat_simplex_3_GL2_16.csv};
        \addlegendentry{$L^2$, $\Delta t=16/N$}
        \addplot[dotted, purple]
        table [x=N,y=hdivsig, col sep=comma]{mixed_heat_simplex_3_GL2_1.csv};
        \addlegendentry{$\hdiv$, $\Delta t=1/N$}
        \addplot[loosely dotted, orange] 
        table [x=N,y=hdivsig, col sep=comma]{mixed_heat_simplex_3_GL2_4.csv};
        \addlegendentry{$\hdiv$, $\Delta t=4/N$}
        \addplot[densely dotted, pink]
        table [x=N,y=hdivsig, col sep=comma]{mixed_heat_simplex_3_GL2_16.csv};
        \addlegendentry{$\hdiv$, $\Delta t=16/N$}
        \addplot [domain=32:128] {0.9/pow(x,2)} node[above, midway, yshift=-1pt, anchor=south west] {$\mathcal{O}(h^2)$};
        \addplot [domain=32:128] {0.008/pow(x,3)} node[above, midway, yshift=-1pt, anchor=south west] {$\mathcal{O}(h^{3})$};
      \end{loglogaxis}
    \end{tikzpicture}
    \label{subfig:gl2m}
    \end{subfigure}
    \hspace{0.04\textwidth}
    \begin{subfigure}[l]{0.475\textwidth}
      \caption{RadauIIA(3)}
      \begin{tikzpicture}[scale=0.88]
      \begin{loglogaxis}[xlabel={$N$}, ylabel={$\|u-u_h\|/\|u\|$},
             ylabel near ticks, ymax=1.e-3, ymin=1.e-7,
             legend pos=south west, legend style={font=\tiny} ,
             cycle list name=color list]
        \addplot[dashed, red]
        table [x=N,y=l2u, col sep=comma]{mixed_heat_simplex_3_R3_1.csv};
        \addlegendentry{$L^2$, $\Delta t=1/N$}
        \addplot[densely dashed, blue] 
        table [x=N,y=l2u, col sep=comma]{mixed_heat_simplex_3_R3_4.csv};
        \addlegendentry{$L^2$, $\Delta t=4/N$}
        \addplot[loosely dashed, green]
        table [x=N,y=l2u, col sep=comma]{mixed_heat_simplex_3_R3_16.csv};
        \addlegendentry{$L^2$, $\Delta t=16/N$}
        \addplot[dotted, purple]
        table [x=N,y=hdivsig, col sep=comma]{mixed_heat_simplex_3_R3_1.csv};
        \addlegendentry{$\hdiv$, $\Delta t=1/N$}
        \addplot[loosely dotted, orange] 
        table [x=N,y=hdivsig, col sep=comma]{mixed_heat_simplex_3_R3_4.csv};
        \addlegendentry{$\hdiv$, $\Delta t=4/N$}
        \addplot[densely dotted, pink]
        table [x=N,y=hdivsig, col sep=comma]{mixed_heat_simplex_3_R3_16.csv};
        \addlegendentry{$\hdiv$, $\Delta t=16/N$}
        \addplot [domain=32:128] {8/pow(x,3)} node[above, midway, yshift=-1pt, anchor=south west] {$\mathcal{O}(h^{3})$};
      \end{loglogaxis}
    \end{tikzpicture}
    \label{subfig:r3m}
    \end{subfigure}
    \label{fig:mixedheat}
    \caption{Relative $L^2$ and $\hdiv$ errors at $t=2$ for the mixed heat equation on an $N \times N$ mesh divided into right triangles using various time stepping schemes and CFL numbers.}
\end{figure}
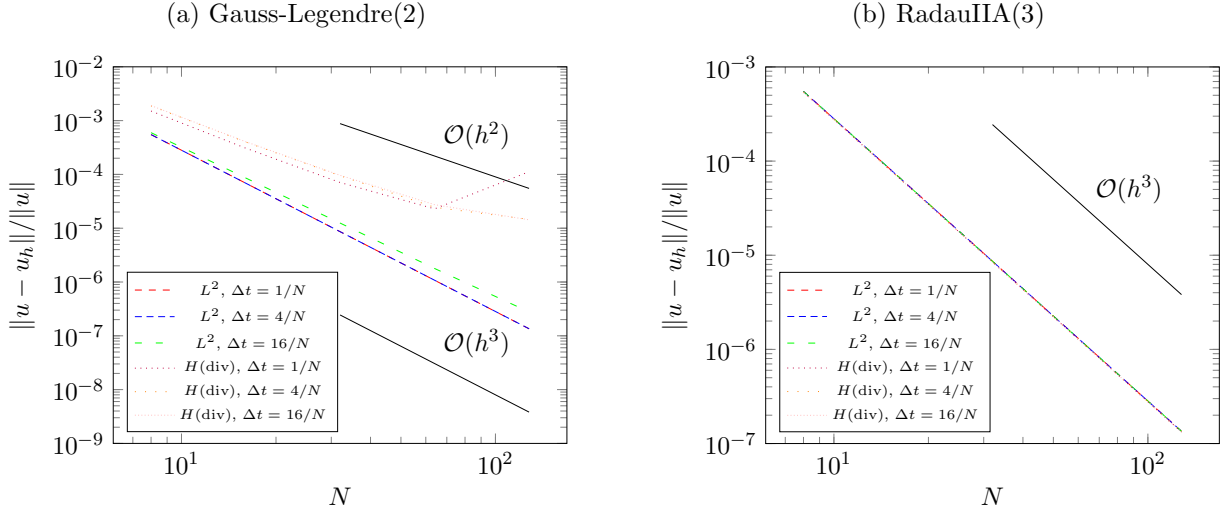

We also consider the mixed form of the heat equation, now dividing the $N \times N$ mesh into right triangles and using third-order Raviart-Thomas elements for $V_h$ with discontinuous polynomials of degree three for $W_h$.  Hence, both variables $u_h$ and $\sigma_h$ should be third-order accurate in space, as should $\nabla \cdot \sigma_h$.  In Figure~\ref{subfig:gl2m}, we see order reduction for a two-stage Gauss-Legendre method, although no order reduction for the three-stage RadauIIA method is observed in Figure~\ref{subfig:r3m}.
From this, we see that the extra formal order of accuracy for Gauss-Lobatto methods does not materialize unless one also takes a very small time step.  We also obtained satisfactory results with somewhat lower accuracy with the two-stage RadauIIA method (but at a lower cost) and the three-stage LobattoIIIC method (at a comparable cost).  Note that the $\hdiv$ error for $\sigma$ and $L^2$ error for $u$ are similar but not identical for the RadauIIA method.  This approximation of the flux variable with comparable accuracy to the potential is seen as a major advantage of mixed methods.

\subsection{Preconditioning}
We consider the implementation of the preconditioner~\eqref{eq:P}.
We see that the diagonal
blocks to be solved are quite similar to a backward Euler method, and
so we choose to solve them inexactly using a single V-cycle of algebraic
multigrid (AMG)~\cite{adams2004}. Suitable solver options for Firedrake, using an additive field split (i.e.~a block diagonal
preconditioner), are given in Listing~\ref{kentgamg}.

%\begin{figure}
%\begin{lstlisting}[captionpos=b, language=Python, caption={Simple PETSc solver options selecting GMRES as the Krylov solver and \texttt{gamg} as the preconditioner.}, label={gamg}]
%  params = {"mat_type": "aij",
%            "snes_type": "ksponly",
%            "ksp_type": "gmres",
%            "pc_type": "gamg"}
%\end{lstlisting}
%\end{figure}
\begin{figure}
\begin{lstlisting}[captionpos=b, language=Python, caption={Possible solver options for a block preconditioner~\eqref{eq:P}.  Here, the inverse of each block diagonal is approximated by a single application of \texttt{gamg} rather than applied via an inner Krylov iteration.}, label={kentgamg}]
    params = {"mat_type": "aij",
              "snes_type": "ksponly",
              "ksp_type": "gmres",
              "ksp_monitor": None,
              "pc_type": "fieldsplit",         
              "pc_fieldsplit_type": "additive"
            }

    per_field = {"ksp_type": "preonly",
                 "pc_type": "gamg"}

    for s in range(butcher_tableau.num_stages):
        params["fieldsplit_%s" % (s,)] = per_field
\end{lstlisting}
\end{figure}
We note that the option for \lstinline{mat_type} was added to
Firedrake in~\cite{kirby2018solver}, and controls whether a nested matrix, a
monolithic sparse matrix, or a ``matrix free'' context is built.  We
also point out that it is straightforward to use different schemes
for each block (for example, when the first stage is explicit, one
only needs a simple method to invert a mass matrix).

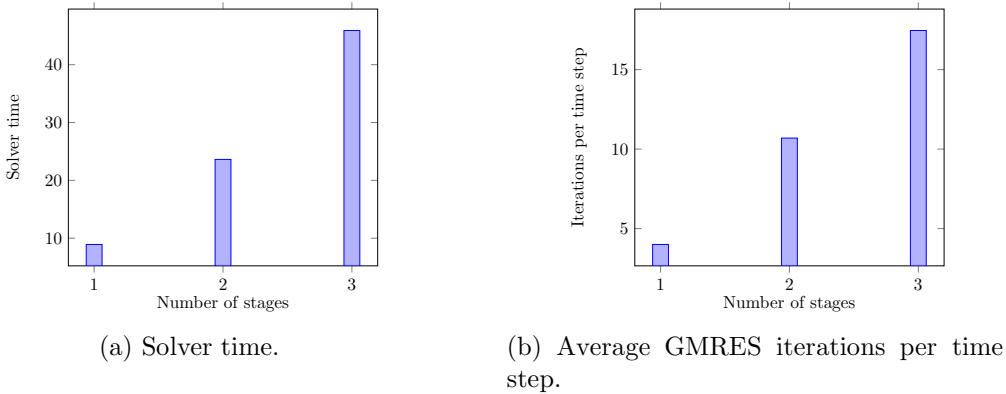
\begin{figure}
  \pgfplotstableread[col sep=comma,]{heat_pc.csv}\datatable
\begin{subfigure}[t]{0.4\textwidth}
\centering
    \begin{tikzpicture}[scale=0.6]
      \begin{axis}[ybar,
        xtick=data, xlabel={Number of stages},
        xticklabels from table={\datatable}{k}, ylabel={Solver time}]
        \addplot table [x=k, y={SNES}]{\datatable};
      \end{axis}
    \end{tikzpicture}
      \caption{Solver time.}
      \label{timeperstage}
\end{subfigure}
\hspace{0.04\textwidth}
\begin{subfigure}[t]{0.4\textwidth}
\centering
    \begin{tikzpicture}[scale=0.6]
      \begin{axis}[ybar,
        xtick=data, xlabel={Number of stages},
        xticklabels from table={\datatable}{k}, ylabel={Iterations per time step}]
        \addplot table [x=k, y={ITS}]{\datatable};
      \end{axis}
    \end{tikzpicture}
      \caption{Average GMRES iterations per time step.}
      \label{its}
\end{subfigure}
\caption{Performance of RadauIIA methods for the heat equation on a $128\times 128$ mesh of $Q^2$ elements.  With $k=1$ we use GMRES preconditioned with a single sweep of \texttt{gamg}, while for $k=2,3$ we use GMRES with an additive field split preconditioner using \texttt{gamg} on each diagonal block.}
\label{fig:timeheatpc}
\end{figure}

For our experiment, we integrate the heat equation from $T=0$ to $T=1$ on a $128 \times 128$ mesh of squares.  We use $Q^2$ elements in space and RadauIIA($k$) methods in time with $k=1,2,3$ and a time step of $\Delta t=0.078125$.  In Figure~\ref{fig:timeheatpc}, we report the total solver time over all time steps and the average GMRES iteration count per time step.  We note that RadauIIA(2) is approximately three times as expensive as RadauIIA(1)/backward Euler.
Each GMRES iteration is about twice as expensive and more iterations are required at each time step.  
Moreover, RadauIIA(3) is observed to be about five times as expensive as backward Euler.  Considering the much greater accuracy obtained by a fifth-order method over a first, we consider this a significant result.

Step-doubling via Richardson extrapolation gives a classical method for building a higher-order method out of a first order method like backward Euler.  To obtain a second-order method in time, one can first advance the solution from $t^n$ to $t^n + \Delta t$ by a single step of size $\Delta t$.  Then, one also takes two time steps of size ${\Delta t}/{2}$ to go from $t^n$ to $t^n + \Delta t$.  A simple linear combination of the two results at $t^n +\Delta t$ gives a second order method.  However, this requires three backward Euler steps to acquire a single second-order time step.  Empirically, our preconditioning strategy requires about the cost of three backward Euler steps for a two-stage method.  For RadauIIA, this gives a method with stiff order two, but full order three for non-stiff problems.

If one wishes to obtain a third-order method by repeated step doubling, one must set up a tableau with four backward Euler steps of size $\tfrac{\Delta t}{4}$, two steps of size $\tfrac{\Delta t}{2}$ and one of size $\Delta t$ -- a total of seven backward Euler steps.  By comparison, our three stage method (empirically) requires a cost comparable to only backward Euler evaluations.  Three-stage RadauIIA has stiff order three and is fifth order for non-stiff problems, so this is a clear win over step doubling.

\subsection{The wave equation}
We now turn to the wave equation~\eqref{eq:mixedwave}.  Here, our experiment fixes a $10 \times 10$ mesh of the unit square subdivided into right triangles with $RT_2 \times DG_1$ spatial discretization  and considers the effect of the time-stepping method and step size on energy conservation.  The semidiscrete system is known~\cite{kirbykieu} to exactly preserve the energy
\begin{equation}
E(t) = \frac{1}{2} \left( \| u \|^2 + \| \sigma \|^2 \right),
\end{equation}
so that $E(t) = E(0)$ for all $t > 0$.  

In particular, we pick $\Delta t = c / N$ for $c=1, 5, 10$ so that $\Delta t=0.1, 0.5, 1.0$.  We consider the two lowest-order methods each of Gauss-Legendre, LobattoIIIC, and RadauIIA.  Our experiment is simple: we integrate from $t=0$ to $t=10$ and compute the ratio of the final energy at $t=10$ to the initial energy.  These ratios are shown in Table~\ref{table:energy}.  Although Gauss-Legendre methods fared poorly relative to L-stable methods for the heat equation, the reverse is true here.  We see the exact energy conservation obtained for the Gauss-Legendre methods, while significant damping occurs for the LobattoIIIC and RadauIIA methods.

\begin{table}
  \begin{center}
  \begin{tabular}{c|ccc} \toprule
     & \multicolumn{3}{c}{$\Delta t$} \\ \midrule
    method & 0.1 & 0.5 & 1.0 \\ \midrule
    Gauss-Legendre(1) & 1 & 1 & 1 \\
    Gauss-Legendre(2) & 1 & 1 & 1 \\ \midrule
    LobattoIIIC(2) & 3.79$\times 10^{-1}$ & 9.75$\times 10^{-18}$ & 1.17$\times 10^{-20}$ \\
    LobattoIIIC(3) & 9.99$\times 10^{-1}$ & 5.19$\times 10^{-2}$ & 3.65$\times 10^{-9}$ \\ \midrule
    RadauIIA(1) & 1.50$\times 10^{-8}$ & 3.40$\times 10^{-16}$ & 6.79$\times 10^{-14}$ \\
    RadauIIA(2) & 9.00$\times 10^{-1}$ & 7.10$\times 10^{-4}$ & 3.77$\times 10^{-7}$\\ \bottomrule
  \end{tabular}
  \end{center}
\caption{Energy conservation for the next-to-lowest order mixed method for the wave equation~\eqref{eq:mixedwave}.  We take a $10\times 10$ mesh divided into right triangles and advance in time to $t=10$ using various time-stepping strategies with various time steps.  We see that the symplectic Gauss-Legendre methods with 1 and 2 stages conserve the energy to machine precision.  However, the other methods do not conserve energy and become very dissipative as the time step is increased.}
\label{table:energy}
\end{table}

So far, all of our examples have used fully implicit methods, but here we point out how we can make use of PETSc options to achieve greater efficiency for DIRKs.  Although this example is in the context of the wave equation, it applies equally to any setting in which a DIRK is appropriate. 

Suppose we use the Qin/Zhang 2-stage DIRK given by Butcher tableau
\begin{equation}
  \begin{array}{c|cc}
    1/4 & 1/4 & 0  \\
    3/4 & 1/2 & 1/4  \\ \hline
    & 1/2 & 1/2
  \end{array}
\end{equation}
which is second-order, A-stable, and symplectic \cite{mei1992diagonally}.

When the Butcher matrix $A$ is lower-triangular, the system matrix~\eqref{eq:blockmatgen} is block lower triangular, and so one may solve for each stage variable in succession via block forward substitution.  This can be implemented as a PETSc preconditioner so that no special care is required in Irksome to handle DIRK methods differently.  One obtains a block lower triangular preconditioner by means of a multiplicative field split~\cite{brown2012composable}, and if that preconditioner is applied exactly to a block lower triangular system, the system will be solved exactly.  

Because Firedrake already has a mixed system ($u$ and $\sigma$), we take care to specify that the system be blocked as a $2\times 2$ system
combining the variables for each stage rather than the default of a $4
\times 4$ split for each separate field.  This can be controlled via PETSc options, and a solver specified for each block as in Listing~\ref{fig:dirkparams}.  While our simple example applies a direct method on each diagonal block, more advanced solvers could also be employed.  Repeating the experiments above for the wave equation with this Butcher tableau and these options gives identical energy behavior to the Gauss-Legendre methods in Table~\ref{table:energy}, although each time step simply consists of solving the diagonal blocks and performing forward substitution (all done inside the PETSc \texttt{FIELDSPLIT} preconditioner).

\begin{figure}
\begin{lstlisting}[captionpos=b, language=Python, caption={Simple DIRK-appropriate parameters for the mixed wave equation allowing one to solve for each stage in succession.  A direct method is used on each stage.}, label={fig:dirkparams}]
  params = {"mat_type": "aij",
            "snes_type": "ksponly",
            "ksp_type": "preonly",
            "pc_type": "fieldsplit",
            "pc_fieldsplit_type": "multiplicative"}
  params["pc_fieldsplit_0_fields"] = "0,1"
  params["pc_fieldsplit_1_fields"] = "2,3"
  per_field = {"ksp_type": "preonly",
               "pc_type": "lu"}
  for i in range(butcher_tableau.num_stages):
      params["fieldsplit_%d" % i] = per_field
\end{lstlisting}
\end{figure}

\subsection{Nonlinear examples}
We conclude our examples with two nonlinear problems.  The
time-dependent Navier-Stokes equations highlight a nonlinear DAE-type
system, and the Benjamin-Bona-Mahony dispersive wave model illustrates
a nonlinear Sobolev-type equation.  Irksome's high-level interface
works unchanged in these examples.

\subsubsection{Navier-Stokes}
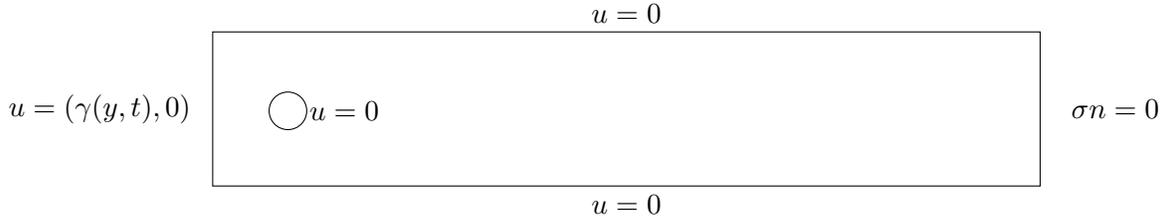
\begin{figure}
  \begin{center}
  \begin{tikzpicture}[scale=5]
    \draw (0,0) rectangle (2.2, 0.41);
    \draw (.2, .2) circle (0.05);
    \node[] at (1.1, 0.46) {$u=0$};
    \node[] at (1.1, -0.05) {$u=0$};
    \node[] at (.35, .2) {$u=0$};
    \node[] at (2.4, 0.205) {$\sigma n = 0$};
    \node[] at (-0.3, 0.205) {$u=\left(\gamma(y, t), 0\right)$};
  \end{tikzpicture}
  \end{center}
  \caption{Computational domain for flow past cylinder, with boundary conditions indicated on each part of the boundary.}
  \label{fig:domain}
\end{figure}

We consider the drag and lift calculations from the two-dimensional benchmark flow around a cylinder proposed in~\cite{john2004reference,schafer1996benchmark}.  The domain is given by $\Omega=[0,2.2] \times [0, 0.41] \backslash B_r(0.2, 0.2)$, with radius $r=0.05$, and is shown in Figure~\ref{fig:domain}. The density is taken as $\rho=1$ and kinematic viscosity is $\nu = 10^{-3}$.  No-slip conditions are imposed on the top and bottom of the pipe and the cylinder.  Natural boundary conditions (no-stress) are imposed on the outflow right end, and a parabolic profile is posed on the inflow boundary on the left end:
\begin{equation}
  u(t, y) = \left( \frac{4 U(t) y(0.41 - y)}{0.41^2}, 0 \right) \equiv \gamma(y, t),
\end{equation}
where $U(t) = 1.5 \sin\left(\tfrac{\pi t}{8}\right)$.  This configuration corresponds to a Reynolds number of 100.  The problem is integrated over one time period, from $t=0$ until $t=8$.

\begin{figure}[h]
  \begin{subfigure}[l]{0.475\textwidth}
    \begin{tikzpicture}[scale=0.88]
      \begin{axis}[xlabel={$t$}, ylabel={Drag}]
        \addplot[dotted, red]
        table [x=T,y=D, col sep=comma]{reference_6.csv.reduced};
        \addlegendentry{Reference}
        \addplot[dashed, green]
        table [x=t,y=CD, col sep=comma]{results_2_0.00125.dat.reduced};
        \addlegendentry{Computed}        
      \end{axis}
    \end{tikzpicture}
    \caption{Drag}
    \label{subfig:drag}
  \end{subfigure}
  \hspace{0.04\textwidth}
  \begin{subfigure}[l]{0.475\textwidth}
    \begin{tikzpicture}[scale=0.88]
      \begin{axis}[xlabel={$t$}, ylabel={Lift},
          legend pos=north west]
        \addplot[dotted, red]
        table [x=T,y=L, col sep=comma]{reference_6.csv};
        \addplot[dashed, green]
        table [x=t,y=CL1, col sep=comma]{results_2_0.00125.dat};
        \legend{Reference, Computed}
      \end{axis}
    \end{tikzpicture}
    \caption{Lift}
    \label{subfig:lift}    
  \end{subfigure}
  \caption{Lift and drag calculations using Firedrake and Irksome compared to benchmark results in~\cite{john2004reference,schafer1996benchmark}.  Excellent agreement is obtained with half as many degrees of freedom and a time-step twice as large, owing to the use of higher-order methods.}
  \label{fig:draglift}
\end{figure}
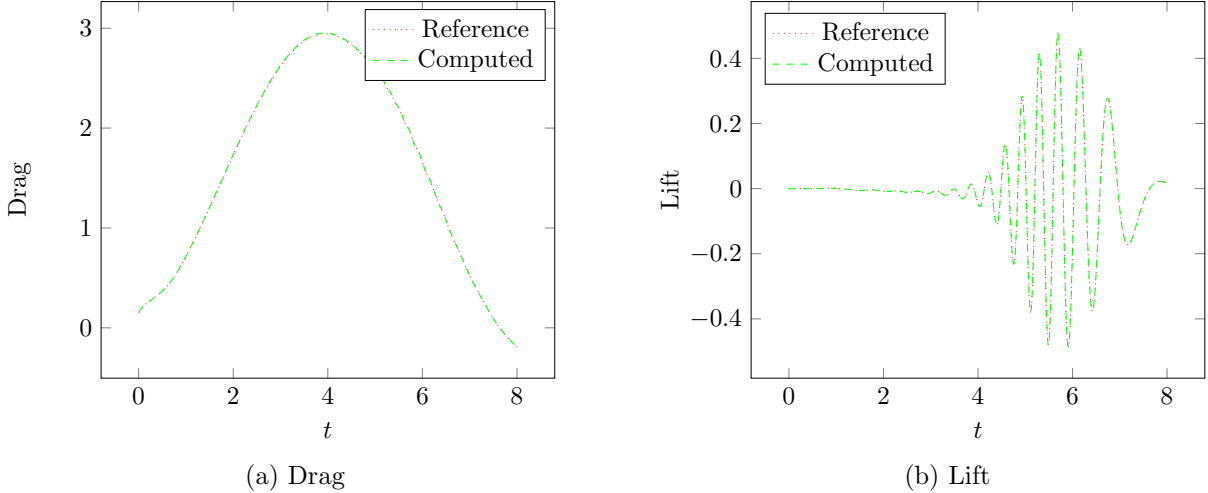

We model this geometry using curvilinear triangles (quadratic mappings capture the geometry exactly).  Using Firedrake's integration with OpenCascade, we are able to mesh the geometry with \texttt{gmsh}~\cite{geuzaine2009gmsh} and refine the mesh in a geometrically conforming way.  Our sample run uses a mesh with 14,640 cells and 7,546 vertices and a Scott-Vogelius discretization~\cite{scott1985norm} with $P^4$ velocities and discontinuous $P^3$ pressures.  This amounts to about 382k global degrees of freedom.  This higher order discretization also suggests a higher order time-discretization,  and we use the two-stage RadauIIA method with $\Delta t= 0.00125$.  
By contrast, the benchmark values to which we compare are obtained on a finer mesh (133,120 straight quadrilateral elements with 133,952 vertices) and lower-order discretization (biquadratic velocity and discontinuous linear pressure) with about 667k degrees of freedom.  They use Crank-Nicolson time-stepping with $\Delta t=1/1600 = 0.000625$, half the size of our time step.  In Figure~\ref{fig:draglift}, we report the computed drag and lift for our method over time compared to the benchmark values, and excellent agreement is attained.

As a remark, Scott-Vogelius discretizations enforce the divergence-free condition pointwise in affine geometry, but this is not the case with our curvilinear mesh.  While the exact theory of this method is still not understood, we have found it to be stable.  The $L^2$ norm of the divergence at each time step was found to range between about $10^{-5}$ and $10^{-3}$, which is comparable to the difference with the reference drag and lift values.  If an exactly divergence-free method is required, one could reduce the geometry representation to affine or perhaps switch to an $\hdiv$/$L^2$ discretization~\cite{cockburn2007note}.

\subsubsection{Benjamin-Bona-Mahony}
We now return to the Benjamin-Bona-Mahony equation~\eqref{eq:bbm}.
With the \lstinline{Dt} operator added from Irksome, we can express the weak form of the semidiscrete equation in UFL by
\begin{lstlisting}
F = (inner(Dt(u), v) * dx + inner(u.dx(0), v) * dx
     + inner(u * u.dx(0), v) * dx + inner((Dt(u)).dx(0), v.dx(0)) * dx)
\end{lstlisting}

We consider this problem on the space interval $[0, 100]$ divided into $N=1000$ intervals (giving mesh size $h=0.1$) and pose periodic boundary conditions.  We take the initial condition as
\[
u(x, 0) = 3 \frac{c^2}{1-c^2} \sech^2 \frac{1}{2} \left( c x + \delta \right),\
\]
where $\delta$ is chosen so that the bump is centered at $x=40$ and $c=\frac{1}{2}$.   This is the initial condition for a right-traveling solitary wave with velocity $\tfrac{1}{1-c^2} = \tfrac{4}{3}$.  We integrate from time 0 to time $T=18$, at which the wave has traveled 24 units to the right.

\begin{figure}[h]
  \begin{subfigure}[l]{0.475\textwidth}
    \begin{tikzpicture}[scale=0.88]
      \begin{axis}[xlabel={$x$}, ylabel={$u(x, t)$},
             legend pos=north west]
        \addplot[dashed]
        table [x=x,y=u, col sep=comma]{bbmGL2_IC.csv};
        \addlegendentry{$u_h(x,0)$}
        \addplot[dotted]
        table [x=x,y=uex, col sep=comma]{bbmGL2_final.csv};
        \addlegendentry{$u_h(x,18)$}
      \end{axis}
    \end{tikzpicture}
    \label{subfig:bbmicfc}
    \caption{Initial condition and numerical solution at $t=18$.}
  \end{subfigure}
  \hspace{0.04\textwidth}
  \begin{subfigure}[l]{0.475\textwidth}
    \begin{tikzpicture}[scale=0.88]
      \begin{axis}[xlabel={$x$}, ylabel={$u(x, 18) - u_h(x, 18)$},
             legend pos=north east]
        \addplot[dashed]
        table [x=x,y=err, col sep=comma]{bbmGL2_final.csv};
      \end{axis}
    \end{tikzpicture}
    \label{subfig:bbmerror}
    \caption{Difference between exact and computed solution at time $t=18$.}
  \end{subfigure}
  \caption{Simulation of the Benjamin-Bona-Mahony problem using $P^1$ elements and GL(2) time-stepping methods.}
  \label{fig:bbmplot}
\end{figure}
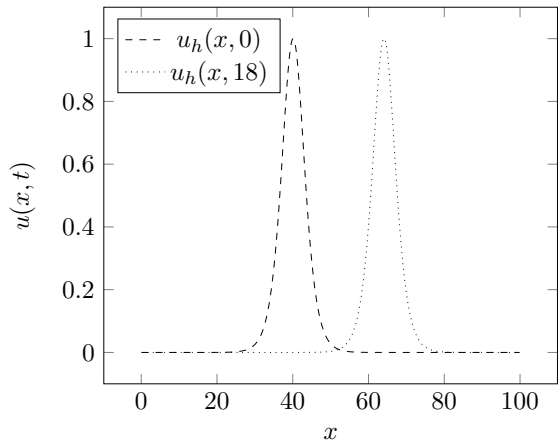
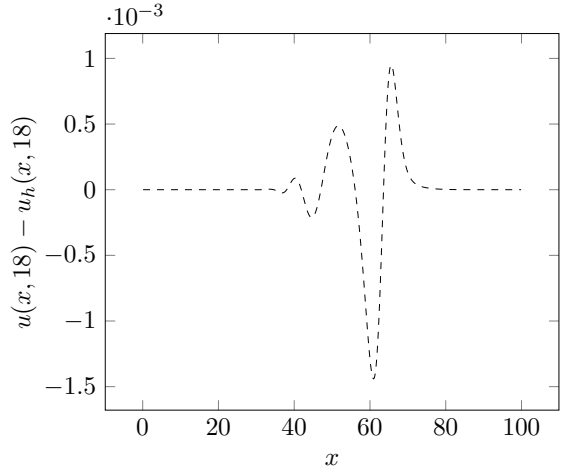

We discretize the problem with a standard Galerkin method using $P^1$ finite elements and, given the Hamiltonian nature of the problem, Gauss-Legendre time stepping methods.  With the implicit midpoint rule, we considered both $\Delta t = h$ and $\Delta t = 10h$.  We found that the former case gave us about $0.15\%$ relative $L^2$ error at $T=18$, but the latter case with a large time step gave us greater than $10\%$ relative error.  Using the fourth-order, two-stage Gauss-Legendre method with $\Delta t = 10h$ actually gave us  $0.14\%$ relative error at the final time -- less than the lower order method with a much smaller time step.  Figure~\ref{fig:bbmplot} shows the solution at the initial and final times as well as the pointwise error between the true and computed solutions at $t=18$.

We also demonstrate conversion of the invariants in Table~\ref{table:bbm_invariant}.  Both $I_1$ and $I_2$ were conserved to $\mathcal{O}(10^{-15})$ for both GL methods.  $I_3$ drifted slightly with both methods, on the order of $10^{-3}$ for the one-stage method (still, less than a 1\% drift) and on the order of $10^{-6}$ for the two-stage method.  While Gauss-Legendre methods are symplectic and preserve quadratic invariants~\cite{sanz1988runge, hairer2006geometric}, very strong conditions are required for them to preserve cubic invariants as well~\cite{iserles2000preserving}.

\begin{table}
  \begin{center}
    \begin{tabular}{c|c}
      \toprule
      GL(1) & GL(2) \\ \midrule
      \begin{tabular}{c|ccc}
        $t$ & $I_1(t)/I_1(0)$ & $I_2(t)/I_2(0)$ & $I_3(t)/I_3(0)$ \\ \midrule
        0 & 1 & 1 & 1 \\
        6 & 1 & 1 & 0.9984516 \\
        12 & 1 & 1 & 0.9963326 \\
        18 & 1 & 1 & 0.9948969
      \end{tabular}
      & 
      \begin{tabular}{c|ccc}
        $t$ & $I_1(t)/I_1(0)$ & $I_2(t)/I_2(0)$ & $I_3(t)/I_3(0)$ \\ \midrule
        0 & 1 & 1 & 1 \\
        6 & 1 & 1 & 0.9999996 \\
        12 & 1 & 1 & 0.9999992 \\
        18 & 1 & 1 & 0.9999990
      \end{tabular} \\
    \end{tabular}
  \end{center}
  \caption{Conservation of invariants for the BBM equation \eqref{eq:bbm} using one (left) and two (right) stage Gauss-Legendre methods.  Invariant $I_1$ is linear and $I_2$ is quadratic, and these are exactly preserved.  The cubic invariant $I_3$ is not exactly conserved, but only very small relative deviation from the initial quantity is observed.}
  \label{table:bbm_invariant}
\end{table}

\section{Conclusions and future work}
\label{sec:conc}

Irksome offers the opportunity for many future directions.  IMEX-type RK methods can be helpful with processes with multiple time scales (e.g.~diffusion and advection) operating on the same variables.  Some additional extensions to UFL, such as the form labeling present in Gusto~\cite{gusto}, would be helpful in allowing users to segregate the terms to be treated implicitly and explicitly. More recently, general linear methods~\cite{butcher2006general} generalize both RK and multi-step methods.  It is conceivable that both multistep and general linear methods could be implemented by transforming UFL in the manner we have proposed here for RK methods.

\bibliographystyle{plain}
\bibliography{references.bib}

\end{document}